\documentclass[12pt]{amsart}

\usepackage{graphicx}
\usepackage{amsthm}
\usepackage{amsmath, amssymb}
\swapnumbers

\theoremstyle{plain}
\newtheorem{Thm}{Theorem}

\newtheorem{Coro}[Thm]{Corollary}
\newtheorem{Lem}[Thm]{Lemma}

\theoremstyle{definition}
\newtheorem{Def}[Thm]{Definition}

\begin{document}

\title[Bounding the stable genera of Heegaard splittings]{Bounding the stable genera of Heegaard splittings from below}

\author{Jesse Johnson}
\address{\hskip-\parindent
        Department of Mathematics \\
        Yale University \\
        PO Box 208283 \\
        New Haven, CT 06520 \\
        USA}
\email{jessee.johnson@yale.edu}


\thanks{Research supported by NSF MSPRF grant 0602368}

\begin{abstract}
We describe for each postive integer $k$ a 3-manifold with Heegaard surfaces of genus $2k$ and $2k-1$ such that any common stabilization of these two surfaces has genus at least $3k-1$.  We also show that for every positive $n$, there is a 3-manifold that has $n$ pairwise non-isotopic Heegaard splittings of the same genus all of which are stabilized.
\end{abstract}

\maketitle

\section{Introduction}

A \textit{Heegaard splitting} for a compact, connected, closed, orientable 3-manifold $M$ is a triple $(\Sigma, H^-, H^+)$ where $\Sigma \subset M$ is a compact, connected, closed, orientable, separating surface and $H^-, H^+ \subset M$ are handlebodies such that $H^- \cup H^+ = M$ and $\partial H^- = \Sigma = \partial H^+ = H^- \cap H^+$.  We will say that two Heegaard splittings are \textit{isotopic} if there is an ambient isotopy taking one of the surfaces to the other.

A \textit{stabilization} of a Heegaard splitting is a new splitting constructed by taking a connect sum of the original splitting with a Heegaard splitting of $S^3$.  Reidemeister~\cite{reid} and Singer~\cite{sing} showed independently that given two Heegaard splittings of the same 3-manifold, there is a third Heegaard splitting, called a \textit{common stabilization}, that is isotopic to a stabilization of each of the initial splittings.

The \textit{stable genus} of two Heegaard splittings is the genus of their smallest common stabilization.  Many examples are known of pairs of Heegaard splittings whose stable genus is $p+1$, where $p$ is the larger of the two initial genera.  It has been a long standing problem to find pairs of Heegaard splittings whose stable genus is higher than this.  We prove the following:

\begin{Thm}
\label{mainthm}
For every $k > 1$, there is a 3-manifold with Heegaard splittings of genus $2k-1$ and $2k$ such that the stable genus of these two Heegaard splittings is $3k-1$.
\end{Thm}

This is proved in Section~\ref{thm1proofsect}.  The examples contain incompressible tori and are therefore not hyperbolic.  However, the construction can be modified to produce atoroidal 3-manifolds with Heegaard splittings of genera $2k$ and $2k-2$ whose stable genus is $3k-2$ for each $k$.  The details of modifying the construction are left to the reader.

Moriah and Sedgwick~\cite{morsedg} have asked whether there is a closed 3-manifold with a weakly reducible Heegaard splitting of non-minimal genus.  In these examples, both Heegaard splittings are weakly reducible.  The genus $2k$ Heegaard splitting has non-minimal genus, so this gives a positive answer to their question.

This paper is a continuation and a generalization of an earlier paper~\cite{me:stabs} by the same author and we will refer to this paper for a number of key Lemmas.  Because this earlier paper deals with a simpler case, the reader may want to review it before reading this paper.  The method of proof in both papers is motivated by Hass, Thompson and Thurston's paper~\cite{htt:stabs}.  They use a hyperbolic geometry argument to show that there exist Heegaard splittings such that the smallest stabilization in which the handlebodies can be interchanged by an isotopy has genus twice that of the original.

David Bachman~\cite{bachman} has recently announced similar examples using different techniques; where we use bicompressible surfaces to compare two Heegaard splittings, he uses incompressible surfaces.

Theorem~\ref{mainthm} implies that the 3-manifold it describes has two stabilized Heegaard splittings that have the same genus but are not isotopic.  By generalizing the construction, we can find 3-manifolds with arbitrarily many stabilized but non-isotopic Heegaard splittings.

\begin{Thm}
\label{mainthm2}
For every $n \geq 1$, there is a 3-manifold with $n$ stabilized Heegaard spittings of the same genus such that no two of them are isotopic.
\end{Thm}

This is proved in Section~\ref{thm2proofsect}.  I would like to thank Andrew Casson, Joel Hass and Abby Thompson for helpful conversations.

\section{Sweep-outs and graphics}
\label{sweepsect}

A \textit{sweep-out} for a compact, orientable 3-manifold $M$ is a smooth function $f : M \rightarrow [-1,1]$ such that each of $f^{-1}(-1)$ and $f^{-1}(1)$ is the union of a graph in $M$ and a collection of boundary components of $M$, while for $t \in (-1,1)$, $f^{-1}(t)$ is a connected, closed surface parallel to $f^{-1}(0)$.  The sets $f^{-1}(-1)$ and $f^{-1}(1)$ are called the \textit{spines} of $f$ and their union contains all of $\partial M$.  

A \textit{stable function} between smooth manifolds $M$ and $N$ is a smooth function $\phi : M \rightarrow N$ such that in the space $C^\infty(M,N)$ of smooth functions from $M$ to $N$, there is a neighborhood $N$ around $\phi$ in which each function is isotopic to $\phi$.  A Morse function is a stable function from a smooth manifold to $\mathbf{R}$ and one can think of stable functions as a generalization of Morse theory to functions whose range has dimension greater than one.

Let $M'$ and $M''$ be 3-dimensional submanifolds of a 3-manifold $M$ and assume $M' \cap M''$ is a non-empty 3-dimensional submanifold $M^*$.  Let $f : M' \rightarrow [-1,1]$ be a sweep-out for $M'$ and $g : M'' \rightarrow [-1,1]$ a sweep-out for $M''$. The product of their restrictions to $M^*$ is a smooth function $f \times g : M^* \rightarrow [-1,1] \times [-1,1]$.  (That is, we define $(f \times g)(x) = (f(x),g(x))$.)  

In the case when $M' = M'' = M$, Kobayashi~\cite{Kob:disc} has shown that after an isotopy of $f$ and $g$, we can assume that $f \times g$ is a stable function on the complement of the spines of $f$ and $g$.  An almost identical argument in the more general case implies that after an isotopy of $f$ and $g$, their product will be stable on the complement in $M^*$ of their spines.

The local behavior of stable functions between dimensions two and three has been classified~\cite{mather} and coincides with the classification by Cerf~\cite{cerf:strat} that was used by Rubinstein and Scharlemann~\cite{rub:compar} to compare Heegaard splittings using pairs of sweep-outs.

\begin{figure}[htb]
  \begin{center}
  \includegraphics[width=4in]{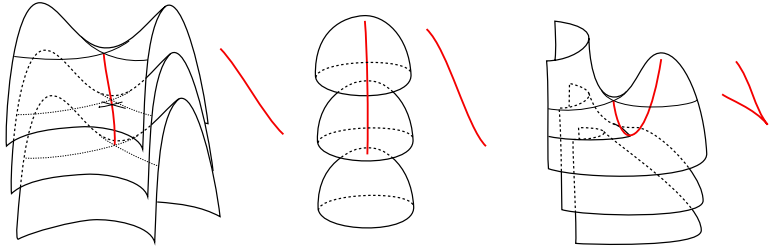}
  \caption{Edges of the graphic are formed by points where the level surfaces of $f$ are tangent to level surfaces of $g$.  Level surfaces of $g$ are shown with $f$ the height function and its level surfaces horizontal planes.}
  \label{grfig}
  \end{center}
\end{figure}

At each point in the complement of the spines, the differential of the map $f \times g$ is a linear map from $\mathbf{R}^3$ to $\mathbf{R}^2$.  This map will have a one dimensional kernel for a generic point in $M$.  The \textit{discriminant set} for $f \times g$ is the set of points where the derivative has a two or three dimensional kernel.  (In fact, all the critical points in a stable function in these dimensions have two dimensional kernels.)  Mather's classification of stable functions~\cite{mather} implies that the discriminant set in this case will be a one dimensional smooth submanifold in the complement in $M$ of the spines.  It consists of all the points where a level surface of $f$ is tangent to a level surface of $g$.  Some examples are shown in Figure~\ref{grfig}.  (For a more detailed description see~\cite{Kob:disc} or~\cite{rub:compar}.)  

The function $f \times g$ sends the discriminant to a graph in $[-1,1] \times [-1,1]$ call the \textit{Rubinstein-Scharlemann graphic} (or just the \textit{graphic} for short).  The parts of the graphic corresponding to the tangencies in Figure~\ref{grfig} are shown next to the surfaces.  The vertices in the interior of the graphic are valence four (crossings) or valence two (cusps).  The vertices in the boundary are valence one or two.

The pre-image in $f \times g$ of an arc $[-1,1] \times \{s\}$ is the level set $g^{-1}(s)$ and the restriction of $f$ to this level surface is a function $\phi_s$ with critical points in the levels where the arc $[-1,1] \times \{s\}$ intersects the graphic as well as possibly at the levels $-1$ and/or $1$.  The same is true if we switch $f$ and $g$.  

\begin{Def}
The function $f \times g$ is \textit{generic} if $f \times g$ is stable on the complement in $M^*$ of the spines of $f$ and $g$ and each arc $\{t\} \times [-1,1]$ or $[-1,1] \times \{s\}$ contains at most one vertex of the graphic.
\end{Def}

If for a given $s \in [-1,1]$ the arc $[-1,1] \times \{s\}$ does not intersect any vertices then every critical point of $\phi_s$ will be non-degenerate and away from $-1$ and $1$ no two critical points will be in the same level.  In other words, $\phi_s$ will be \textit{Morse away from} $-1$ and $1$.  If the arc passes through a vertex then in the levels other than $-1$ and $1$, $\phi_s$ will either have a degenerate critical point or two non-degenerate critical points at the same level.  We will say that such a $\phi_s$ is \textit{near-Morse} away from $-1$ and $1$.

\section{Labeling the graphic}
\label{labelsect}

\begin{Def}
A \textit{compression body} is a connected 3-manifold homeomorphic to a regular neighborhood $H$ of the union of a connected graph $K$ properly embedded in a 3-manifold $M$ and every component of $\partial M$ that contains a vertex of $K$.  The union of $K$ and the boundary components is called a \textit{spine} for $H$.  
\end{Def}

Note that a handlebody is a compression body formed from a graph that is disjoint from $\partial M$.  We will write $\partial_- H = \partial H \cap \partial M$ and $\partial_+ H = \partial H \setminus \partial_- H$.  Note that $\partial_+ H$ is connected and has higher genus than every component of $\partial_- H$.  For a 3-manifold $M$ with boundary, we define a Heegaard splitting of $M$ to be a triple $(\Sigma, H^-, H^+)$ where $\Sigma \subset M$ is a closed surface and $H^-, H^+$ are compression bodies with $\partial_+ H^- = \Sigma = \partial_+ H^+ = H^- \cap H^+$ and $H^- \cup H^+ = M$.  For such a Heegaard splitting, we have $\partial M = \partial_- H^- \cup \partial_+ H^+$.

Let $M$ be a compact, connected, closed, orientable 3-manifold and $S \subset M$ a connected, closed, two-sided, separating surface in $M$.  (This argument can be generalized to most non-separating surfaces, but for simplicity we will restrict our attention to separating surfaces.)  Throughout the paper, we will assume that $S$ has genus at least two, $M$ is irreducible and $S$ is not contained in a ball in $M$.  

Let $N$ be a closed regular neighborhood of $S$.  Because $S$ is 2-sided, $N$ is homeomorphic to $S \times [-1,1]$ and we will identify $S \times [-1,1]$ with $N$ such that $S \times \{0\} = S$.  Define $N^- = S \times [-1,0]$ and $N^+ = S \times [0,1]$.

\begin{Def}
The surface $S$ is bicompressible if there are embedded disks $D^-, D^+ \subset M$ such that $D^- \cap S = \partial D^-$ and $D^+ \cap S = \partial D^+$ are essential loops in $S$, $D^- \cap N$ is contained in $N^-$ and $D^+ \cap N$ is contained in $N^+$.
\end{Def}

In other words, $S$ is bicompressible if there are essential compressing disks on both sides of the surface.  If $S$ is bicompressible then both boundary components of $N$ are compressible into $M \setminus N$.  Let $D$ be a compressing disk such that $\partial D$ is an essential loop in $\partial N$ and the interior of $D$ is contained in $M \setminus N$.  Let $N'$ be the union of $N$ and a closed regular neighborhood of $D$.  This set is the result of compressing $\partial N$ into $M \setminus N$.  

If there is a compressing disk for $\partial N'$ then we can repeat the process.  Let $M_S$ be the result of compressing the boundary of $N$ maximally into $M$, i.e. repeating this process until $\partial M_S$ cannot be compressed further into $M \setminus M_S$.  (The surface $\partial M_S$ will, however, be compressible into $M_S$.)  If any component of $\partial M_S$ is a sphere then it bounds a ball disjoint from $S$ in $M$ (because $M$ is irreducible and $S$ is not contained in a ball) and we will add all such balls into $M_S$.

By construction, the surface $S$ will be separating in $M_S$.  Let $M_S^-$ and $M_S^+$ be the closures of the components of $M_S \setminus S$ such that $N^- \subset M_S^-$ and $N^+ \subset M_S^+$.  Each of these sets is a compression body so $S$ is a Heegaard surface for $M_S$.  In order to understand $M_S$, we need the following Lemma about compression bodies:

\begin{Lem}
\label{cbodyincomplem}
If $H$ is a compression body then $\partial_- H$ is incompressible in $H$.  Conversely, every closed, positive genus, incompressible surface embedded in $H$ is parallel to a component of $\partial_- H$.
\end{Lem}

This lemma can be proved by noting that $\partial_+ H$ can be compressed until the resulting 3-manifold is homeomorphic to $\partial_- H \times [0,1]$.  The boundary components of this manifold are incompressible.  Any incompressible surface in $H$ and any compressing disk for $\partial_- H$ can be isotoped disjoint from the compression that produced $\partial_- H \times [0,1]$.  The details of this proof are left to the reader.

\begin{Lem}
\label{msuniquelem}
The set $M_S$ is uniquely defined up to isotopy (i.e. it doesn't matter what disks you compress along, as long as they're maximal).
\end{Lem}

\begin{proof}
Let $M^1_S$ be the submanifold resulting from compressing along one set $\mathbf{D}^1$ of disks and let $M^2_S$ be the result of compressing along a second set, $\mathbf{D}^2$.  Isotope $\mathbf{D}^1$ to intersect $\partial M^2_S$ minimally.  If an innermost loop of intersection in $D \subset \mathbf{D}^1$ is essential in $\partial M^2_S$ then $\partial M^2_S$ can be compressed further, contradicting the maximality assumption on $\mathbf{D}^2$.  Thus an innermost disk intersects $\partial M^2_S$ in a trivial loop.  Since $M$ is irreducible, this disk can be isotoped into $M^2_S$.  

After isotoping $\mathbf{D}^1$ to minimize the intersection, the disks will be contained in $M^2_S$.  Thus we can isotope $M^1_S$ into $M^2_S$.  If $\partial M^1_S$ is compressible in $M^2_S \setminus M^1_S$ then $\mathbf{D}^1$ is not maximal.  The set $M^1_S$ is a union of compression bodies $M^{1-}_S$, $M^{1+}_S$ along their positive boundaries.  The set $M^2_S$ is also a union of compression bodies $M^{2-}_S$, $M^{2+}_S$.  By Lemma~\ref{cbodyincomplem}, $\partial M^{1-}_S$ is not compressible into $M^{1-}_S$, so $\partial_- M^{1-}_S$ is incompressible in $M^{2-}_S$.  Each component of $\partial M^{1-}_S$ is contained in $M^{2-}_S$.  So by Lemma~\ref{cbodyincomplem}, each component of $\partial_- M^{1-}_S$ must be parallel to a component of $\partial_- M^{2-}_S$.  

The set $\partial M^{1-}_S$ is separating in $M^{2-}_S$.  No proper subset of $\partial M^{2-}_S$ is separating, so there is at least one component of $\partial M^{1-}_S$ parallel to each component of $M^{2-}_S$.  Conversely, two surfaces parallel to any component of $M^{2-}_S$ bound a surface cross interval component.  Because $M^{1-}_S$ has no such component, there must be exactly one component of $\partial M^{1-}_S$ parallel to each component of $\partial M^2_S$.  The same argument applies to $M^{1+}_S$ and $M^{2+}_S$.  We can thus isotope $M^1_S$ so that it is contained in $M^2_S$ and their boundaries coincide, implying $M^1_S$ is isotopic to $M^2_S$.
\end{proof}

We will say that a sweep-out $f : M_S \rightarrow [-1,1]$ \textit{represents} $(S, M_S^-, M_S^+)$ if $f$ can be isotoped so that $f^{-1}(-1)$ is contained in $M_S^-$, $f^{-1}(1)$ is contained in $M_S^+$ and $f^{-1}(t)$ is parallel in $M_S$ to $S$ for each $t \in (-1,1)$.  

Let $(\Sigma, H^-, H^+)$ be a Heegaard splitting for $M$.  We will say that a sweep-out $g : M \rightarrow [-1,1]$ \textit{represents} $(\Sigma, H^-, H^+)$ if $g$ can be isotoped so that $g^{-1}(-1)$ is contained in $H^-$, $g^{-1}(1)$ is contained in $H^+$ and $g^{-1}(s)$ is parallel to $\Sigma$ for each $s \in (-1,1)$.  This is a special case of the discussion in Section~\ref{sweepsect} with $M_\Sigma = M$.  The intersection of $M$ and $M_S$ is precisely $M_S$ so as in the Section~\ref{sweepsect}, we can isotope $f$ and $g$ so that $f \times g$ is generic.

If $S$ is in fact a Heegaard surface then $M_S = M$ and we are in the situation considered in~\cite{me:stabs}.  The reader can check that the definitions below are simply generalizations of the same terms defined in the previous paper.

Given sweep-outs $f$ for $(S, M^-_S, M^+_S)$ and $g$ for $(\Sigma, H^-, H^+)$ as above we will define $S_t = f^{-1}(t)$, $\Sigma_s = g^{-1}(s)$, $H^-_s = g^{-1}([-1,s])$ and $H^+_s = g^{-1}([s,1])$.

\begin{Def}
For $t, s \in [-1,1]$, we will say that $S_t$ is \textit{mostly above} $\Sigma_s$ if each component of $S_t \cap H^-_s$ is contained in a disk in $S_t$.  Similarly, $S_t$ is \textit{mostly below} $\Sigma_s$ if each component of $S_t \cap H^+_s$ is contained in a disk in $S_t$.
\end{Def}

Define $R_a \subset (-1,1) \times (-1,1)$ as the set of all ordered pairs $(t,s)$ such that $S_t$ is mostly above $\Sigma_s$.  Similarly, define $R_b \subset (-1,1) \times (-1,1)$ as the set of all ordered pairs $(t,s)$ such that $S_t$ is mostly below $\Sigma_s$.  Note that the upper boundary of $R_a$ and the lower boundary of $R_b$ are contained in the graphic.

\begin{Def}
We will say that $g$ \textit{spans} $f$ if $f \times g$ is generic and for some values $s, t_-, t_+ \in (-1,1)$, $S_{t_-}$ is mostly below $\Sigma_s$, while $S_{t_+}$ is mostly above $\Sigma_s$.  Moreover, we will say that $g$ spans $f$ \textit{positively} if $t_- < t_+$, or \textit{negatively} if $t_- > t_+$.
\end{Def}

Note that $g$ will span $f$ if and only if there is a horizontal arc $[-1,1] \times \{s\}$ that intersects both $R_a$ and $R_b$.  If $f \times g$ is generic and there is no such arc, then we have our next definition.

\begin{Def}
The sweep-out $g$ \textit{splits} $f$ if $f \times g$ is generic and there is a horizontal arc $[-1,1] \times \{s\}$ that is disjoint from both $R_a$ and $R_b$.
\end{Def}

\begin{figure}[htb]
  \begin{center}
  \includegraphics[width=2.5in]{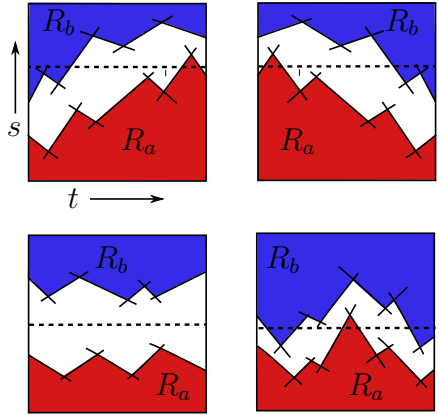}
  \put(-150,145){$R_b$}
  \put(-20,145){$R_b$}
  \put(-130,55){$R_b$}
  \put(-65,55){$R_b$}
  \put(-120,100){$R_a$}
  \put(-60,100){$R_a$}
  \put(-107,7){$R_a$}
  \put(-34,7){$R_a$}
  \put(-163,105){$s$}
  \put(-140,78){$t$}
  \caption{Clockwise from the top left, the graphics correspond to pairs of sweep-outs $f$, $g$ such that (1) $g$ spans $f$ positively, (2) $g$ spans $f$ negatively, (3) $g$ both spans $f$ positively and spans $f$ negatively and (4) $g$ splits $f$.  The dotted line represents the arc $[-1,1] \times \{s\}$.}
  \label{spanningfig}
  \end{center}
\end{figure}

These conditions are illustrated in Figure~\ref{spanningfig}.   

We can extend them directly to a pair of bicompressible surfaces as follows:

\begin{Def}
We will say that $(\Sigma, H^-, H^+)$ \textit{spans} $(S, M^-_S, M_S^+)$ positively or negatively if $(\Sigma, H_-, H_+)$ and $(S, M^-_S, M_S^+)$ are represented by sweep-outs $g$ and $f$ such that $g$ spans $f$ positively or negatively, respectively.  We will say that $(\Sigma, H^-, H^+)$ \textit{splits} $(S, M^-_S, M_S^+)$ if $(\Sigma, H_-, H_+)$ and $(S, M^-_S, M_S^+)$ are represented by sweep-outs $g$ and $f$ such that $g$ splits $f$. 
\end{Def}

As an example, we note the following Lemma, which is a generalization of Lemma~11 in~\cite{me:stabs}.  The proof is identical to the proof of Lemma~11 in~\cite{me:stabs} and will be left as an exercise for the reader.

\begin{Lem}
\label{stabspanslem}
If a Heegaard splitting $(\Sigma, H^-, H^+)$ spans a bicompressible surface $(S, M_S^-, M_S^+)$ positively then every stabilization of $(\Sigma, H^-, H^+)$ spans $(S, M_S^-, M_S^+)$ positively.  If $(\Sigma, H^-, H^+)$ spans $(S, M_S^-, M_S^+)$ negatively then every stabilization of $(\Sigma, H^-, H^+)$ spans $(S, M_S^-, M_S^+)$ negatively.
\end{Lem}

\section{Amalgamating Heegaard splittings}
\label{amalgsect}

Consider 3-manifolds $M_1$ and $M_2$ and Heegaard splittings $(\Sigma_1, H_1^-, H_1^+)$, $(\Sigma_2, H_2^-, H_2^+)$ for $M_1$, $M_2$, respectively such that $\partial_- H_1^+$ is a non-empty subset of $\partial M_1$ (possibly all of $\partial M_1$) homeomorphic to $\partial_- H_2^+ \subset M_2$.  Let $M$ be the union of $M_1$ and $M_2$ with $\partial_- H_1^+$ glued to $\partial_- H_2^+$ by some homeomorphism.  The images in $M$ of $H_1^-$, $H_1^+$, $H_2^-$ and $H_2^+$ form what's called a \textit{generalized Heegaard splitting} with $\Sigma_1$ and $\Sigma_2$ the \textit{thick surfaces} and the image of $\partial_- H^-_1$ the \textit{thin surface}, as in~\cite{sch:thin}.  Let $F \subset M$ be image of $\partial_- H^-_1$.

We will construct a Heegaard splitting $(\Sigma_3, H^-_3, H^+_3)$ for $M$ from this generalized Heegaard splitting as follows:  Because $H^+_1$ is a compression body, it can be decomposed into a submanifold $A_1$ homeomorphic to $\partial_- H^+_1 \times [0,1]$ and a submanifold $B_1$ homeomorphic to a collection of balls (1-handles) such that $A_1 \cup B_1 = H^+_1$ and $A_1 \cap B_1$ is a collection of disks.  Let $A_2, B_2$ be a similar collection for $H^+_2$.  

Let $C_1 \subset A_1 (= \partial_- H_1^+ \times [0,1])$ be the submanifold $(A_1 \cap B_1) \times [0,1]$.  This is a collection of balls, each of which intersects $F = \partial_- H^+_1$ in a single disk.  Let $C_2 = (A_2 \cap B_2) \times [0,1]$ be equivalent set in $B_2$.  Because $C_1$ and $C_2$ intersect $F$ in two collections of disks, we can choose them (by changing the identification between $A_i$ and $\partial_- H_i^+ \times [0,1]$) so that they are disjoint.

Note that $A_1 \cup C_1$ and $A_2 \cup C_2$ are collections of balls such that each collection intersects $F$ in a union of disjoint disks in their boundaries.  The set $B_2 \setminus C_2$ is a regular neighborhood of a subsurface of $\partial_+ H^-_2$ so the closure $X^+$ of $H^-_2 \cup (B_2 \setminus C_2)$ is a compression body such that $\partial_+ X^+$ contains $F \setminus C_2$.  In particular, $\partial_+ X^+$ contains $(A_1 \cup C_1) \cap F$ so $X^+ \cup (A_1 \cup C_1)$ is a compression body.  Define $H_3^+ = X^+ \cup (A_1 \cup C_1)$.  Similarly, we can define a second compression body $H_3^- = X^- \cup (A_2 \cup C_2)$ where $X^-$ is the closure of $H^-_1 \cup (B_1 \setminus C_1)$.  These compression bodies coincide along their positive boundaries, so they determine a Heegaard splitting for $M$.

\begin{Def}
The Heegaard splitting constructed above is called an \textit{amalgamation} of $(\Sigma_1, H_1^-, H^+_1)$ and $(\Sigma_2, H_2^-, H^+_2)$ along the surface $F$.
\end{Def}

Note that we have labeled the amalgamation so that the first handlebody $H^-_1$ of $(\Sigma_1, H^-_1, H^+_1)$ is contained in the first handlebody $H^-_3$ of $(\Sigma_3, H^-_3, H^+_3)$  but the first handlebody $H^-_2$ of $(\Sigma_2, H^-_2, H^+_2)$ is contained in the second handlebody $H^+_3$ of $(\Sigma_3, H^-_3, H^+_3)$.  We will say that this ordering of the compression bodies \textit{agrees with} $(\Sigma_1,H^-_1, H^+_1)$.  The ordering $(\Sigma_3, \Sigma_3^+, \Sigma_3^-)$ will agree with $(\Sigma_2, H^-_2, H^+_2)$ (as well as with $(\Sigma_1, H^+_1, H^-_1)$).

The genus of the amalgamation can be calculated directly.  The reader can check that the genus of $\Sigma_3$ is equal to the sum of the genera of $\Sigma_1$ and $\Sigma_2$ minus the sum of the genera of the components of $F$.

\begin{Lem}
\label{amalgspanlem}
Let $M$ be a 3-manifold, $(S,M^-_S, M^+_S)$ a bicompressible surface and $(\Sigma, H^-, H^+)$ a Heegaard splitting for $M$ such that $(\Sigma, H^-, H^+)$ spans $(S, M_S^-, M_S^+)$ positively.  Let $M'$ be the result of gluing a 3-manifold $M''$ to $M$ along $\partial_- H^+$.  Let $(\Sigma', H'^-, H'^+)$ be an amalgamation of $(\Sigma, H^-, H^+)$ with a Heegaard splitting for $M''$.  If the labeling of $(\Sigma', H'_1, H'_2)$ agrees with $(\Sigma, H^-, H^+)$ then $(\Sigma', H'^-, H'^+)$ will span $(S, M_S^-, M_S^+)$ positively.
\end{Lem}

\begin{proof}
Let $f$ be a sweep-out for $(S, M_S^-, M_S^+)$ and $g$ a sweep-out for $(\Sigma, H^-, H^+)$ such that $g$ spans $f$ positively.  Let $s$, $t_-, t_+$ be values such that $S_{t_-}$ is mostly below $\Sigma_s$ while $S_{t_+}$ is mostly above $\Sigma_s$.

To form an amalgamation, we construct the compression body $H'^-$ for the amalgamation by isotoping $H^-$ and then taking its union with a collection of balls that intersect $\partial_+ H^-$ in disks (1-handles).  This collection of balls is isotopic to a regular neighborhood of a graph with one or more vertices in $\partial_+ H^-$.  Thus we can isotope $H'^-$ so that it is the union of $H^-_s$ with a regular neighborhood of a graph.  

If we choose this regular neighborhood small enough then the intersection of $S_{t_-}$ with $H'^-$ will differ from its intersection with $H^-_s$ by a collection of disks.  Thus $S_{t_-}$ will be mostly below $\Sigma'$.  Similarly, $S_{t_+}$ will be mostly above this isotoped $\Sigma'$.  If we choose a sweep-out for $\Sigma'$ such that $\Sigma'_0$ is this isotoped surface then this sweep-out will span $f$ positively.  Thus $(\Sigma', H'^-, H'^+)$ spans $(S, M_S^-, M_S^+)$ positively.
\end{proof}

Consider a connected, closed, orientable surface $F$.  The 3-manifold $F \times [0,1]$ has two boundary components so there are two types of Heegaard splittings of $F \times [0,1]$, determined by whether the two boundary components are on opposite sides or the same side of the Heegaard surface.  The surface $F \times \{\frac{1}{2}\}$ is a Heegaard surface that splits the boundary components and we will call this a \textit{trivial} Heegaard splitting of $F \times [0,1]$.

Consider an arc $\alpha = \{x\} \times [0,1]$ for some point $x \in F$.  The closure $H_1^+$ of a regular neighborhood of $\alpha \cup (F \times \{0,1\})$ is a compression body and the closure $H_1^-$ of its complement is a handlebody.  If $\Sigma_1$ is their common boundary then $(\Sigma_1, H^-_1, H^+_1)$ is a Heegaard splitting for $F \times [0,1]$ and we will call this a \textit{minimal non-trival} Heegaard splitting.

Let $M$ be a 3-manifold with a single boundary component $F$.  Let $(\Sigma, H^-, H^+)$ be a Heegaard splitting for $M$ such that $H^-$ is a handlebody and $\partial_- H^+ = F$.  If we glue $F \times [0,1]$ to $M$ along $\partial M$ and one of the components of $\partial (F \times [0,1])$ then the resulting 3-manifold is homeomorphic to $M$.  We can thus produce a new Heegaard splitting for $M$ by amalgamating $\Sigma$ with a Heegaard splitting for $F \times [0,1]$.

\begin{Def}
The amalgamation of $(\Sigma, H^-, H^+)$ with a minimal non-trivial Heegaard splitting for $F \times [0,1]$ is called a \textit{boundary stabilization} of $\Sigma$ along $F$.
\end{Def}

The observant reader may note that we gave a different definition of boundary stabilization in~\cite{me:stabs}.  These definitions are in fact equivalent (i.e. the constructions produce isotopic Heegaard splittings.)  For a 3-manifold with one boundary component, we would like to keep the convention that $H^-$ is a handlebody.  Thus we will label the handlebodies of a boundary stabilization to agree with the Heegaard splitting of $F \times [0,1]$ rather than with the original Heegaard splitting.  Thus Lemma~\ref{amalgspanlem} above implies the following, which is a generalization of Lemma~13 in~\cite{me:stabs} to this situation:

\begin{Lem}
\label{bstablspanslem}
Let $M$ be a 3-manifold, $S$ a bicompressible surface and $(\Sigma, H^-, H^+)$ a Heegaard splitting for $M$ such that $(\Sigma, H^-, H^+)$ spans $(S, G^-, G^+)$ positively.  Then a boundary stabilization of $(\Sigma, H^-, H^+)$ will span $(S, G^-, G^+)$ negatively.
\end{Lem}

\section{Spanning sweep-outs}

In this section we will prove the following generalization of Lemma~14 in~\cite{me:stabs}.  Recall that we have assumed $M$ is irreducible.

\begin{Lem}
\label{bothspanlem}
If $\Sigma$ spans $S$ then $\Sigma$ is an amalgamation along $S$.  If there is a sweep-out $g$ representing $\Sigma$ and a sweep-out $f$ representing $S$ such that $g$ spans $f$ both positively and negatively then $\Sigma$ is an amalgamation along two copies of $S$, such that the Heegaard surface for the $S \times [0,1]$ component does not separate the two copies of $S$.
\end{Lem}

The following Lemma is precisely Lemma~15 in~\cite{me:stabs}. We will refer the reader to that paper for the proof.

\begin{Lem}[Lemma 15 in \cite{me:stabs}]
\label{compressthenamalglem}
If there is a sequence of compressions that turn a Heegaard surface $\Sigma$ into a surface $F$ then $\Sigma$ is an amalgamation along $F$.
\end{Lem}

We will prove Lemma~\ref{bothspanlem} by combining Lemma~\ref{compressthenamalglem} with the following:

\begin{Lem}
\label{spanthencompresslem}
If $\Sigma$ spans $S$ then there is a sequence of compressions of $\Sigma$ that turn $\Sigma$ into a surface with a component isotopic to $S$.  If a sweep-out representing $\Sigma$ spans $S$ positively and negatively then $\Sigma$ compresses down to a surface containing at least two components isotopic to $S$.
\end{Lem}

\begin{proof}
We will prove the second case of the lemma, when one sweep-out spans the other with both signs.  The first case follows from a very similar, but even simpler, argument.  The proof follows the arguments in Section~5 of~\cite{me:stabs} and we will refer to a number of Lemmas from that section in this proof.

Let $f$ be a sweep-out representing $(S, M^-_S, M^+_S)$ and $g$ a sweep-out representing $(\Sigma, H^-, H^+)$.  Assume $f \times g$ is generic and that $g$ spans $f$ both positively and negatively.  By Lemma~17 of~\cite{me:stabs}, there is a value $s$ and values $t_- < t_0 < t_+$ such that either $S_{t_-}$ and $S_{t_+}$ are mostly below $\Sigma_s$ and $S_{t_0}$ is mostly above $\Sigma_s$ or $S_{t_-}$ and $S_{t_+}$ are mostly above $\Sigma_s$ and $S_{t_0}$ is mostly below $\Sigma_s$.  (This should appear obvious from the bottom right graphic in Figure~\ref{spanningfig}.)  Without loss of generality, we will assume $S_{t_-}$ and $S_{t_+}$ are mostly below $\Sigma_s$ while $S_{t_0}$ is mostly above $\Sigma_s$.

Define $F_0 = \Sigma_s$.  By definition, each component of $F_0 \cap S_{t_-} = \Sigma_s \cap S_{t_-}$ is a trivial loop in $S_{t_-}$.  Let $\ell \subset \Sigma_s \cap S_{t_-}$ be an innermost loop, bounding a disk in $S_{t_-}$.  If $\ell$ bounds a disk in $F_0$ then the union of the two disks is a sphere, which bounds a ball since $M$ is irreducible.  Let $F_1$ be the result of isotoping $F_0$ across this ball so as to eliminate the component $\ell$ of the intersection.  Otherwise, if $\ell$ is essential in $F_0$ then let $F_1$ be the result of compressing $F_0$ across the disk in $S_{t_-}$.  This also eliminates the loop $\ell$ from the intersection.

Continue the compression process with $S_{t_-}$, then with $S_{t_0}$ and $S_{t_+}$ to form a surface $F_n$ disjoint from $S_{t_-}$, $S_{t_0}$ and $S_{t_+}$.  After each compression, we can decompose the complement of $F_i$ into two sets $G_i^-$, $G_i^+$ such that for each component $F'_i$ of $F_i$, one of the two components of $M \setminus F_i$ adjacent to $F'_i$ is in $G_i^-$ and the other is in $G_i^+$.  Choose the labels so that either $G^-_i \subset G^-_{i-1}$ or $G^+_i \subset G^+_{i-1}$ (depending on which way $F_{i-1}$ is compressed.

Because the intersection of $S_{t_-}$ with $H^+$ is trivial in $S_{t_-}$, the intersection of $S_{t_-}$ with each $G_i^-$ is trivial.  The final $F_n = \partial G_n^+$ is disjoint from $S_{t_-}$ so $S_{t_-}$ must be contained in $G_n^-$.  Similarly, $S_{t_0}$ is contained in $G_n^+$ so $F_n$ separates $S_{t_-}$ from $S_{t_0}$.  Finally, $S_{t_+}$ is contained in $G^-_n$ so $F_n$ separates $S_{t_0}$ from $S_{t_+}$.  

By Lemma~16 of~\cite{me:stabs} (with the roles of $F$ and $S$ reversed), the components of $F_n$ contained in $f^{-1}([t_-,t_0])$ can be compressed further to a surface containing a component isotopic to some $S_t$.  (This is a simple corollary of the classification of incompressible surfaces in $S \times [0,1]$.)  The same is true for the components of $F_n$ in $f^{-1}([t_0,t_+])$.  Thus $F_n$ compresses further to a surface $F_m$ containing two components isotopic to $S$.
\end{proof}

\begin{proof}[Proof of Lemma~\ref{bothspanlem}]
First assume $\Sigma$ spans $S$.  By Lemma~\ref{spanthencompresslem}, $\Sigma$ compresses down to a surface $F$ containing a component isotopic to $S$.  By Lemma~\ref{compressthenamalglem}, this implies that $\Sigma$ is an amalgamation of a generalized Heegaard splitting along $F$.  Let $S' \subset F$ be a component of $F$ isotopic to $S$.  Since $S'$ is separating, the generalized Heegaard splitting determines a generalized Heegaard splitting for each component of $M \setminus S'$.  Amalgamating each of these generalized splittings produces a Heegaard splitting for each component of $M \setminus S$ and $\Sigma$ is an amalgamation along $S'$ of these two Heegaard splittings.

In the case when a sweep-out for $\Sigma$ spans a sweep-out for $S$ in both directions, Lemma~\ref{spanthencompresslem} implies $\Sigma$ compresses down to a surface $F_m$ containing a subsurface $S'$ consisting of two components, each isotopic to $S$.  If $F_m$ contains more than two components isotopic to $S$ then assume $S'$ consists of adjacent components.  This surface $S'$ is separating, so applying the above argument implies $\Sigma$ is an amalgamation along $S'$.  We can reconstruct the Heegaard surface for the $S \times [0,1]$ component of $M \setminus S'$ by pushing $S \times \{0\}$ and $S \times \{1\}$ into the interior of $S \times [0,1]$ and attaching tubes.  At least one tube passes between the two compoents so the resulting surface does not separate the two components of $S'$.
\end{proof}

\section{Splitting sweep-outs}

Given a compact, connected, closed, orientable surface $S$, the \textit{curve complex} $C(S)$ is the simplicial complex whose vertices are isotopy classes of essential simple closed curves in $S$ and whose simplices are pairwise disjoint sets of loops.  The distance $d(\ell^-,\ell^+)$ between simple closed curves $\ell^-,\ell^+$ in $S$ is defined as the length of the shortest edge path in $C(S)$ between the vertices that represent them.  

Given a connected, separating, bicompressible surface $S \subset M$ and a Heegaard splitting $(S, M_S^-, M_S^+)$ for $M_S$, we will define $\mathcal{M}_S^- \subset C(S)$ as the set of all essential simple closed curves that bound disks in $M_S^-$.  Similarly, $\mathcal{M}_S^+ \subset C(S)$ will be the set of all essential simple closed curves that bound disks in $M_S^+$.  The \textit{(Hempel) distance} $d(S)$ is the distance from $\mathcal{M}_S^-$ to $\mathcal{M}_S^+$, i.e. the minimum of $d(\ell^-, \ell^+)$ over all pairs of essential loops such that $\ell^-$ bounds a disk in $M_S^-$ and $\ell^+$ bounds a disk in $M_S^+$.  By Lemma~\ref{msuniquelem}, $M_S^-$ and $M_S^+$ are uniquely determined (up to isotopy) by $S$ so $d(S)$ is uniquely determined by $S \subset M$.

We will prove the following Lemma, which is a generalization of Lemma~19 in~\cite{me:stabs}.  The proof is identical to the more restricted case, but for the sake of completeness, we will give the complete argument in this more general situation.

\begin{Lem}
\label{splitlem}
Let $k$ be the genus of $\Sigma$.  If $\Sigma$ splits $S$ then $d(S) \leq 2k$.
\end{Lem}

This will be proved at the end of the Section after we have established a number of intermediate lemmas.

By definition, if $\Sigma$ splits $S$ then there are sweep-outs $f$ and $g$ representing $S$ and $\Sigma$, respectively, such that $f \times g$ is generic and $g$ splits $f$.  As noted above, the boundaries of the closures of $R_a$ and $R_b$ are edges of the graphic for $f \times g$.  As pointed out in~\cite{inflects}, a horizontal tangency in the graphic for $f \times g$ corresponds to a critical point in the function $g$.  Since $g$ is a sweep-out, it has no critical points away from its spines, so there can be no horizontal tangencies in the interior of the graphic.  Thus the maxima of the upper boundary of $\bar R_a$ and minima of the lower boundary of $\bar R_b$ are vertices of the graphic.  

Let $C$ be the complement in $\{0\} \times (-1,1)$ of the projections of $R_a$ and $R_b$.  This is a (possibly empty) closed interval.  Because $f \times g$ is generic, if $C$ is a single point, $C = \{s\}$, then the arc $[-1,1] \times \{s\}$ must pass through a single vertex of the graphic that is a maximum of $\bar R_a$ and a mimimum of $\bar R_b$.  Let $(t, s)$ be the coordinates of this vertex.  

For arbitrarily small $\epsilon$, the restriction of $g$ to $S_{t+\epsilon}$ is a Morse function.  Moreover, there are two consecutive critical points in the restriction such that each component of the subsurface below any level set below the first saddle is contained in a disk while each component of any subsurface above a level set above the second saddle is contained in a disk.  This is only possible in a torus.  

Since we assumed $S$ has genus at least two, this is a contradiction and the set $C$ must have more than one point.  Since there are finitely many vertices in the graphic and $C$ is a non-trivial interval, there is an $s \in C$ such that the arc $[-1,1] \times \{s\}$ does not pass through a vertex of the graphic.  

\begin{Lem}
\label{nicesmorselem}
If $g$ splits $f$ then there is an $s$ such that $[-1,1] \times \{s\}$ is disjoint from $R_a$ and $R_b$ and the restriction of $f$ to $\Sigma_s \cap M_S$ is a Morse function away from $0$ and $1$ such that each level set in $\Sigma_s$ contains a loop that is essential in the corresponding level set of $f$.
\end{Lem}

\begin{proof}
As above, we can choose $s$ such that $[-1,1] \times \{s\}$ is disjoint from the vertices of the graphic and from $R_a$ and $R_b$.  The restriction of $f$ to $\Sigma_s \cap M_S$ is Morse away from $0$ and $1$ because $[-1,1] \times \{s\}$ does not pass through any vertices of the graphic.  Each level set of the restriction is a collection of level sets in some $S_t$ that bound the intersection of $S_t$ with $H^-_s$ and with $H^+_s$.  Since $S_t$ is neither mostly above nor mostly below $\Sigma_s$, these loops cannot all be trivial in $S_t$.  Thus the level set contains a loop that is essential in $S_t$.
\end{proof}

To simplify the notation, we will assume (by isotoping if necessary) that $\Sigma = \Sigma_s$ for this value of $s$.

If $d(S) \leq 2$ then Lemma~\ref{splitlem} follows immediately, since we assumed $\Sigma$ has genus at least 2.  Thus we will assume $d(S) > 2$.  Bachman and Schleimer~\cite[Claims 6.3 and 6.7]{bsc:bndls} showed that in this case, there is a non-trivial interval $[a,b] \subset [-1,1]$ such that for $t \in [a,b]$, every loop of $S_t \cap \Sigma_s$ that is trivial in $\Sigma_s$ is trivial $S_t$ while for $t < a$, some loop of $S_t \cap \Sigma_s$ bounds an essential disk in $H^-_t$ and for $t > b$, some loop of $S_t \cap \Sigma_s$ bounds an essential disk in $H^+_t$.

Let $a'$ be a regular level of $f|_{\Sigma}$ just above $a$ and let $b'$ be a regular level just below $b$.  Since $a'$ is in the interval $(a,b)$, every component of $\Sigma \cap S_{a'}$ that is trivial in $\Sigma$ is trivial in $S_{a'}$.  The same is true for $S_{b'}$.

An innermost such loop in $\Sigma$ bounds a disk disjoint from $S_{a'}$ and a second disk in $S_{a'}$.  By assumption, $M$ is irreducible so the two disks cobound a ball.  Isotoping the disk in $\Sigma$ across this ball removes the trivial intersection.  By repeating this process with respect to $S_{a'}$ and $S_{b'}$, we can produce a surface $\Sigma'$ isotopic to $\Sigma$ such that each loop $\Sigma' \cap S_{a'}$ and $\Sigma' \cap S_{b'}$ is essential in $\Sigma'$.  Note that this does not change the property that each regular level set of $f|_{\Sigma'}$ contains a loop that is essential in $S_t$.

Let $F$ be the intersection of $\Sigma'$ with $f^{-1}([a',b'])$.  Consider a projection map $\pi$ from $f^{-1}([a',b'])$ onto $S_0$.  The image of a level loop of $f|_{F}$ under $\pi$ is a simple closed curve in $\Sigma_0$.  (Its isotopy class is well defined, even though its image depends on the choice of projection.)

\begin{Lem}
\label{parallelparallellem}
If two level loops of $f|_{F}$ are isotopic in $F$ then their projections are isotopic in $S_0$.
\end{Lem}

\begin{proof}
Any two level loops are disjoint in $F$ so if two level loops are isotopic then they bound an annulus $A \subset F$.  The projection of $A$ into $S_0$ determines a homotopy from one boundary of the image of $A$ to the other.  Thus the projections of the two loops are homotopic in $S_0$.  Homotopic simple closed curves in surfaces are isotopic so the two projections are in fact isotopic.
\end{proof}

Let $L$ be the set of all isotopy classes of level loops of $f|_F$.  These loops determine a pair-of-pants decomposition for $F$.  We will define a map $\pi_*$ from $L$ to the disjoint union $C(S_0) \cup \{0\}$ as follows:  A representative of a loop $\ell \in L$ projects to a simple closed curve in $S_0$.  If the projection is essential then we define $\pi_*(\ell)$ to be the corresponding vertex of $C(S)$.  If the projection is trivial then we define $\pi_*(\ell) = 0$.  By Lemma~\ref{parallelparallellem}, $\pi_*$ is well defined.  

\begin{Lem}
\label{pantsdsjtlem}
If $\ell$ and $\ell'$ are cuffs of the same pair of pants in the complement $F \setminus L$ then their images in $S_0$ are isotopic to disjoint loops.
\end{Lem}

\begin{proof}
Let $\ell, \ell', \ell'' \in L$ be three loops bounding a pair of pants in $F \setminus L$.  There is a saddle singularity in $f|_{F}$ contained in a level component $E$ (a graph with one vertex and two edges) such that $\ell$, $\ell'$ and $\ell''$ are isotopic to the boundary loops of a regular neighborhood of $E$.

The projection of $E$ into $S_0$ is a graph $\pi(E)$ with one vertex and two edges.  The projections of the level loops near $E$ define a homotopy from the projections of representatives of $\ell$, $\ell'$, $\ell''$ into $\pi(E)$.  Since these representatives are simple in $S_0$, they must be isotopic to the boundary components of a regular neighborhood of $\pi(E)$.  Thus $\pi_*(\ell)$ is disjoint from $\pi_*(\ell')$.
\end{proof}

Thus if $\ell$ and $\ell'$ are cuffs of the same pair of pants and their projections are essential in $S_0$ then $\pi_*(\ell)$ and $\pi_*(\ell')$ are connected by an edge in $C(S_0)$.  Define $L' = \pi_*(L) \cap C(S)$.

\begin{Lem}
\label{lconnectedlem}
The set $L'$ is connected and has diameter at most $2k - 2$.
\end{Lem}

\begin{proof}
For each regular value $t \in (a,b)$ of $f|_F$, let $L_t \subset L$ be the set of loops with representatives in $(f|_F)^{-1}(t)$.  The loops in $L_t$ are pairwise disjoint so their projections in $S_0$ are pairwise disjoint.  Moreover, the projection $\pi(L_t)$ contains at least one essential loop, so $L'_t = \pi_*(L_t) \cap C(S)$ is a non-empty simplex in $C(S)$.  If there are no critical points of $f|_F$ between $t$ and $t'$ then the level sets are isotopic, so $L_t = L_{t'}$ and $L'_t = L'_{t'}$.

If there is a single critical point in $f|_F$ between $t$ and $t'$ then $L_t$ may be different from $L_{t'}$.  If the critical point is a central singularity (a maximum or a minimum) then the difference between the level sets is a trivial loop in $F$, so $L'_t = L'_{t'}$.  If the critical point is a saddle then one or two loops in $L_t$ are replaced by one or two loops in $L_{t'}$.  The two or three loops involved cobound the same pair of pants so by Lemma~\ref{pantsdsjtlem} the corresponding loops in $L'_t$ and $L'_{t'}$ are disjoint.  Thus for any values $t, t' \in [a',b']$, there is a path in $L'$ from any vertex of $L'_t$ to any vertex in $L'_{t'}$.  Since $L'$ is the union of all the sets $\{L'_t | t \in [a',b']\}$, $L'$ is connected.

Consider loops $\ell, \ell' \in L$ whose projections are essential in $S_0$.  Since $L'$ is connected, there is a path $\pi_*(\ell) = v_0, v_1,\dots,v_n = \pi_*(\ell')$ in $L' \subset C(S)$.  Assume we have chosen the shortest such path.  Each $v_i$ is the projection of a loop $\ell_i \in L$.  If $\ell_i$ and $\ell_j$ are cuffs of the same pair of pants in $F \setminus L$ then $v_i$ and $v_j$ are distance one in $C(S_0)$.  Since the path is minimal, $i$ and $j$ must be consecutive.  Each pair of pants contains at most two loops in the path.  Each loop in the path, except possibly the first and last loop, is contained in two pairs of pant.  Thus the number of loops is at most one more than the number of pairs of pants in $F \setminus L$.

The number of pairs of pants is at most the negative Euler characteristic of $F$.  Since $\partial F$ is essential in $\Sigma'$, the Euler characteristic of $F$ is greater (less negative) than or equal to that of $\Sigma'$.  The Euler characteristic of $\Sigma'$ is $2 - 2k$ so the path from $\pi_*(\ell)$ to $\pi_*(\ell')$ has length at most $2k - 2$.
\end{proof}

\begin{proof}[Proof of Lemma~\ref{splitlem}]
Assume $\Sigma$ splits $S$.  Let $[a, b] \subset [-1,1]$ be the largest interval such that for $t \in [a,b]$, every loop of $S_t \cap \Sigma_s$ that is trivial in $\Sigma_s$ is trivial $S_t$.   Let $a',b' \in [-1,1]$ be just inside $[a,b]$ as defined above.  Isotope $\Sigma$, as described, to a surface $\Sigma'$ such that $F = \Sigma' \cap f^{-1}([a',b'])$ has essential boundary in $\Sigma'$ and each level set $(f|_F)^{-1}(t)$ contains an essential loop in $S_t$ for $t \in [a',b']$.  

For small enough $t$, the level loops of $f|_{\Sigma'}$ bound disks in $\Sigma'$ or are parallel to loops in $\partial M_S$.  There will be an essential loop bounding a disk if and only if $a > 0$.  In this case, the value $a$ is a critical level of $f|_{\Sigma'}$ containing a saddle singularity.  As above, the projections of the level loops before and after this essential saddle are pairwise disjoint.  By the definition of $a$, the projection of the level loops before the saddle contain a vertex of $\mathcal{M}_S^-$.  The projection of the level set after $a$ is contained in $L'$ so $d(\mathcal{M}_S^-,L') = 1$.  

If $a = 0$ then every essential loop in $L_{a'}$ is parallel to $\partial M_S$.  Because $M^-_S$ is a compression body, there is an essential, properly embedded disk in the compression body $M^-_S$ disjoint from any such loop so again, $d(\mathcal{H}^-,L') = 1$.

A parallel argument implies $d(\mathcal{M}_S^+,L') = 1$.  By Lemma~\ref{lconnectedlem}, the set $L'$ of projections of level loops into $\Sigma_0$ is connected and has diameter at most $2k-2$.  Thus $d(\Sigma) \leq 2k$.
\end{proof}

\section{Isotopies of sweep-outs}

We have assumed $S$ is a separating surface in $M$.  The closure of each component of $M \setminus S$ is a submanifold of $M$ and we can consider the Heegaard genus of each component.  Let $k'$ be the sum of the Heegaard genera.  In this section, we prove the following:

\begin{Lem}
\label{bothspanboundlem}
If $(\Sigma, H^-, H^+)$ spans $(S, M_S^-, M_S^+)$ both positively and negatively then $k \geq \min\{k', \frac{1}{2}d(S)\}$.
\end{Lem}

If $(\Sigma, H^-, H^+)$ spans $(S, M_S^-, M_S^+)$ both positively and negatively then it will be represented by one sweep-out that spans a sweep-out for $S$ positively and another that spans a sweep-out for $S$ negatively.  These sweep-outs will be isotopic (possibly after some handle slides that do not affect the spanning condition) and we would like to understand how the graphic changes during this isotopy.

\begin{Lem}
\label{isotopesweepslem}
Let $g$ and $g'$ be sweep-outs such that $f \times g$ and $f \times g'$ are generic and $g$ is isotopic to $g'$.  Then there is a family of sweep-outs $\{g_r | r \in [0,1]\}$ such that $g = g_0$, $g' = g_1$ and for all but finitely many $r \in [0,1]$, $f \times g_r$ is generic.  At the finitely many non-generic points, there are at most two valence two or four vertices at the same level, or one valence six vertex.
\end{Lem}

The analogous Lemma for isotopies of Morse functions is Lemma 9 in~\cite{me:t3} and Lemma~\ref{isotopesweepslem} can be proved by a similar argument.  We will allow the reader to work out the details.

We will now prove Lemma~\ref{bothspanboundlem}.  The proof is almost identical to that of Lemma~19 in~\cite{me:stabs}, but we will repeat it here to verify that it works in the more general context.  

\begin{proof}[Proof of Lemma~\ref{bothspanboundlem}]
Since $(\Sigma, H^-, H^+)$ spans $(S, M_S^-, M_S^+)$ both positively and negatively, there are sweep-outs $f$, $g$ representing $(S, M_S^-, M_S^+)$ and $(\Sigma, H^-, H^+)$, respectively, such that $g$ spans $f$ positively, as well as sweep-outs $f'$, $g'$ representing the two surfaces such that $g'$ spans $f'$ negatively.  

The sweep-outs $f$ and $f'$ represent the same bicompressible surface $S$.  By Lemma~\ref{msuniquelem}, $M_S$ is unique up to isotopy so we can compose $f'$ with an isotopy of $M$ after which $f$ and $f'$ will have the same domain.  The sweep-outs then represent isotopic Heegaard splittings of the same submanifold of $M$ so there is a sequence of handle slides after which there is an isotopy taking $f'$ to $f$.  The handle slides can be done in an arbitrarily small neighborhood of the original spine so that before the isotopy, $g'$ still spans $f'$ negatively.  By composing $g'$ with this isotopy, we can assume $g'$ spans $f$ negatively.  Because $g$ and $g'$ represent the same Heegaard splitting, they will be isotopic after an appropriate sequence of handle slides that again do not change the fact that $g$ spans $f$ positively and $g'$ spans $f$ negatively.

Consider a continuous family of sweep-outs $\{g_r |  r \in [0,1], g_r \in C^{\infty}(M,\mathbf{R})\}$ such that $g_0 = g$, $g_1 = g'$ and $f \times g_r$ is generic for all but finitely many $r$, as in Lemma~\ref{isotopesweepslem}.  For a generic $r$, $g_r$ either spans $f$ or splits $f$.  If $g_r$ splits $f$ then by Lemma~\ref{splitlem}, $k \geq \frac{1}{2} d(\Sigma)$.    

If $g_r$ spans $f$ with both signs then by Lemma~\ref{bothspanlem}, $\Sigma$ is an amalgamation of a generalized Heegaard splitting along two copies of $S$.  These two copies of $S$ cut $M$ into a component homeomorphic to $S \times [0,1]$ and two components whose Heegaard genera sum to $k'$ (by definition of $k'$).  By  Lemma~\ref{bothspanlem}, the Heegaard surface in the $S \times [0,1]$ component does not separate the two copies of $S$ so it has genus at least twice that of $S$. Thus by the formula described in Section~\ref{amalgsect}, the genus of the amalgamated Heegaard splitting is at least $k'$.  

Thus if $g_r$ splits $f$ or spans $f$ with both signs then $k \geq \min\{\frac{1}{2}d(\Sigma),k'\}$.  We will therefore assume for contradiction that away from the finitely many non-generic values, $g_r$ spans $f$ positively or negatively, but not both.

Since $g_0$ spans $f$ positively and $g_1$ spans $f$ negatively, there must be some non-generic value $r_0$ such that for small $\epsilon > 0$, $g_{r_0-\epsilon}$ spans $f$ positively, while $g_{r_0+\epsilon}$ spans $f$ negatively.  For every small $\epsilon > 0$, the closures of the projections of $R_a$ and $R_b$ at time $r_0 - \epsilon$ intersect in an interval $I^-_\epsilon$.  Since the projections are disjoint at time $r_0$, the limit of the closures of these intervals must contain a single point $s^-$.  Thus the graphic at time $r_0$ must have two vertices at the same level, one of which is a maximum for the upper boundary of $R_a$ and the other a minimum for the lower boundary of $R_b$, as in the middle graphic shown in Figure~\ref{changefig}.
\begin{figure}[htb]
  \begin{center}
  \includegraphics[width=3.5in]{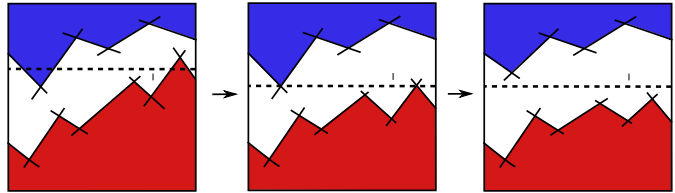}
  \caption{When the graphic goes from spanning positively to not spanning positively, there are two vertices at the same level.}
  \label{changefig}
  \end{center}
\end{figure}

If the vertices in the upper boundary of $R_a$ and the lower boundary of $R_b$ coincide, then this vertex cannot be valence four, as explained above, since $\Sigma$ is not a torus.  The same argument implies that this cannot happen at a valence six vertex either.  Since $g_{r_0 - \epsilon}$ spans $f$ positively, the $s$ coordinate of the vertex in the boundary of $R_a$ must be strictly lower than that the vertex in the boundary of $R_b$.  However, an analogous argument for the graphics at times $r_0 + \epsilon$ implies that the $s$ coordinate of the vertex in the boundary of $R_a$ must be strictly greater than that of the vertex in the boundary of $R_b$.  Since there are at most two vertices at level $s$, this contradiction completes the proof.
\end{proof}

\section{The Proof of Theorem~\ref{mainthm}}
\label{thm1proofsect}

Fix $k \geq 2$ and let $(S_1, G^-_1, G^+_1)$ be a genus $k$ Heegaard splitting of a 3-manifold $M_1$ such that $\partial M_1$ is a torus.  By Hempel's construction~\cite{Hempel:complex}, we can choose $S_1$ so that $d(S_1) > 6k$.  Assume we have done so and let $(S_2, G^-_2, G^+_2)$ be a second genus $k$ Heegaard splitting of a 3-manifold $M_2$ such that $\partial M_2$ is a torus and $d(S_2) > 6k$.  We will follow the convention that $G^-_1$ and $G^-_2$ are handlebodies while $G^+_1$ and $G^+_2$ are compression bodies.

Let $M$ be the union of $M_1$ and $M_2$ identified along their boundaries by some homeomorphism.  By Casson and Gordon's Lemma~\cite[Lemma 1.1]{cass:red}, if $\partial M_1$ is compressible into $M$ then $S_1$ is weakly reducible (i.e. $d(S_1) \leq 1$).  The same holds for $M_2$, $S_2$ so the surface $F = \partial M_1 = \partial M_2$ is incompressible in $M$.  The images in $M$ of $S_1$ and $S_2$ are bicompressible and for $i = 1,2$, we can construct $M_{S_i}$ and a Heegaard splitting $(S_i, M_{S_i}^-, M_{S_i}^+)$ for $M_{S_i}$.  

\begin{Lem}
The set $M_{S_1}$ is isotopic to the image in $M$ of $M_1$.  Similarly, $M_{S_2}$ is isotopic of the image of $M_2$.  Moreover, the distances of $S_1$ and $S_2$ as subsurfaces of $M$ are equal to their distances as Heegaard surfaces in $M_1$, $M_2$, respectively.
\end{Lem}

\begin{proof}
We can compress $S$ in both directions to fill in $M_1$ so we can choose $M_{S_1}$ so that it contains $M_1$.  The boundary of $M_1$ is the torus $F = \partial M_2$.  As noted above, $F$ is incompressible into $M_2 = M \setminus M_1$ so this set of compressions is maximal.  By Lemma~\ref{msuniquelem}, any $M_{S_1}$ that we construct by different compressions is isotopic to $M_1$.  A similar argument for $S_2$ implies $M_{S_2}$ is isotopic to $M_2$.  A loop in $S_i$ will bound a disk in $M_i^\pm$ if and only if it bounds a disk in $M_{S_i}^\pm$, so the distances of $S_1$ and $S_2$ are the same in $M$ as in $M_1$, $M_2$.
\end{proof}

Let $(\Sigma_1, H^-_1, H^+_1)$ be the Heegaard splitting that results from amalgamating the Heegaard surfaces $S_1$ and $S_2$ for $M_1$, $M_2$ along $F$.  Let $(\Sigma_2, H^-_2, H^+_2)$ be the Heegaard splitting that results from boundary stabilizing $S_1$ along $\partial M_1$, then amalgamating the boundary stabilization with $S_2$ along $F$.  (This is the same as boundary stabilizing $S_2$, then amalgamating with $S_1$.)  In both cases, we will order the handlebodies to agree with $(S_2, G^-_2, G^+_2)$.  We will show that the smallest common stabilization of these two Heegaard splittings has genus at least $3k-1$.

\begin{Lem}
The Heegaard splitting $(\Sigma_1, H^-_1, H^+_1)$ spans $(S_1, G^-_1, G^+_1)$ negatively and $(S_2, G^-_2, G^+_2)$ positively.
\end{Lem}

\begin{proof}
The order of the compression bodies of was chosen to agree with $(S_2, G^-_2, G^+_2)$ so by Lemma~\ref{amalgspanlem}, $(\Sigma_1, H^-_1, H^+_1)$ spans $(S_2, G^-_2, G^+_2)$ positively.   If we reverse the order of $H^-_1$ and $H^+_1$ then the order agrees with $S_1$.  Thus $(\Sigma_1, H^+_1, H^-_1)$ spans $(S_1, G^-_1, G^+_1)$ positively.  By the definition of spanning, this implies that $(\Sigma_1, H^-_1, H^+_1)$ spans $(S_1, G^-_1, G^+_1)$ negatively.
\end{proof}

\begin{Lem}
The Heegaard splitting $(\Sigma_2, H^-_2, H^+_2)$ spans both $(S_1, G^-_1, G^+_1)$ and $(S_2, G^-_2, G^+_2)$ positively.
\end{Lem}

\begin{proof}
Again, the ordering of $(\Sigma_2, H^-_2, H^+_2)$ was chosen to agree with $(S_2, G^-_2, G^+_2)$ so by Lemma~\ref{amalgspanlem}, $(\Sigma_2, H^-_2, H^+_2)$ spans  $(S_2, G^-_2, G^+_2)$ positively.  By Lemma~\ref{bstablspanslem}, a boundary stabilization of $(S_1, G^-_1, G^+_1)$ will span $(S_1, G^-_1, G^+_1)$ negatively.  Because the ordering of $(\Sigma_2, H^-_2, H^+_2)$ was chosen to agree with $S_2$, it disagrees with $(S_1, G^-_1, G^+_1)$.  The boundary compression spans $(S_1, G^-_1, G^+_1)$ negatively, so Lemma~\ref{amalgspanlem} implies that $(\Sigma_2, H^-_2, H^+_2)$ spans $(S_1, G^-_1, G^+_1)$ positively.
\end{proof}

We would like to calculate the Heegaard genera of the components of $M \setminus S_1$ and $M \setminus S_2$.  These manifolds have compressible boundary, so we will need the following:

\begin{Lem}
\label{compressgenuslem}
Let $M$ be a 3-manifold with compressible boundary and let $M'$ be the result of removing from $M$ a regular neighborhood of a non-separating, properly embedded, essential disk $D$.  Then the Heegaard genus of $M$ is equal to the Heegaard genus of $M'$ plus one.
\end{Lem}

\begin{proof}
Let $(\Sigma, H^-, H^+)$ be a minimal genus Heegaard splitting of $M$ and assume $\partial D$ is contained in $\partial_- H^+$.  By Casson and Gordon's Lemma~\cite{cass:red}, we can isotope $D$ so that $D \cap \Sigma$ is a single loop.  Let $\Sigma'$ be the image in $M'$ of the result of compressing $\Sigma$ across $D \cap H^+$.  The image of this surface in $M'$ is a Heegaard surface so the Heegaard genus of $M$ is no less than the Heegaard genus of $M'$ minus one.  The reverse of this construction implies the reverse inequality.
\end{proof}

\begin{Coro}
\label{hgencoro}
The Heegaard genus of $M_1 \cup G^+_2$ is equal to $2k-1$, as is the Heegaard genus of $M_2 \cup G^+_1$.
\end{Coro}

\begin{proof}
The boundary of $M_1 \cup G^+_2$ is the surface $\partial_+ G^+_2$ and this surface is compressible into $G^+_2$.  In fact, there are $k-1$ disjoint compressing disks whose union is non-separating and such that compressing along these disks turns $M_1 \cup G^+_2$ into a 3-manifold homeomorphic to $M_1$.
The 3-manifold $M_1$ has a genus $k$ Heegaard surface with distance at least $6k$ so by Scharlemann and Tomova's Thoerem~\cite{st:dist} (See also Corollary~20 of~\cite{me:stabs}) the Heegaard genus of $M_1$ is exactly $k$.  Since we can produce $M_1$ from $M_1 \cup G^+_2$ by sequentially compressing along $k - 1$ non-separating compressing disks, Lemma~\ref{compressgenuslem} implies that the Heegaard genus of $M_1 \cup G^+_2$ is $2k-1$.  An identical proof implies the same for $M_2 \cup G^+_1$.
\end{proof}

Each $S_i$ cuts $M$ into a genus $k$ handlebody (with Heegaard genus $k$) and a second manifold with genus $2k-1$ so the value $k'$ in Lemma~\ref{bothspanboundlem} corresponding to each $S_i$ is $3k-1$.  This is precisely what we need to prove Theorem~\ref{mainthm}.

\begin{proof}[Proof of Theorem~\ref{mainthm}]
Let $(\Sigma'_1, H'^-_1, H'^+_1)$ and $(\Sigma'_2, H'^-_2, H'^+_2)$ be stabilizations of the Heegaard splittings $(\Sigma_1, H^-_1, H^+_1)$ and $(\Sigma_2, H^-_2, H^+_2)$ constructed at the beginning of the section, such that $\Sigma'_1$ is isotopic to $\Sigma'_2$.  By Lemma~\ref{stabspanslem}, $(\Sigma'_1, H'^-_1, H'^+_1)$ spans $(S_1, G^-_1, G^+_1)$ and $(S_2, G^-_2, G^+_2)$ with the same signs as $(\Sigma_1,H^-_1, H^-_1)$.  Thus it spans $(S_1, G^-_1, G^+_1)$ with negatively and spans $(S_1, G^-_1, G^+_1)$ positively, while $(\Sigma'_1, H'^-_1, H'^+_1)$ spans both bicompressible surfaces positively.

The isotopy that takes $\Sigma'_1$ to $\Sigma'_2$ takes $H'^-_1$ to either $H'^-_2$ or $H'^+_2$ and takes $H'^+_1$ to the other handlebody.  If the isotopy takes $H'^-_1$ to $H'^-_2$ and $H'^+_1$ to $H'^+_1$ then $(\Sigma'_1, H'^-_1, H'^+_1)$ spans $(S_1, G^-_1, G^+_2)$ both positively and negatively.  By Lemma~\ref{bothspanboundlem}, this implies that the genus of $\Sigma'_1$ is at least $3k-1$.  Similarly, if the isotopy takes $H'^-_1$ to $H'^+_2$ and $H'^+_1$ to $H'^-_1$ then $(\Sigma'_1, H'^-_1, H'^+_1)$ spans $(S_2, G^-_2, G^+_2)$ both positively and negatively and again we conclude that the genus of $\Sigma'_1$ is at least $3k-1$.
\end{proof}

\section{The Proof of Theorem~\ref{mainthm2}}
\label{thm2proofsect}
In order to find 3-manifolds with more than two non-isotopic stabilized Heegaard splittings of the same genus, we will iterate the construction in the previous section.

Given $n \geq 2$, let $M_1,\dots,M_n$ be 3-manifolds such that each of $M_1$ and $M_n$ has a single torus boundary component while each $M_i$ with $i \neq 1,n$ has two boundary components.  Assume each $M_i$ has a genus four Heegaard surface $S_i$ with distance strictly greater than $6n+8$.  For each $M_i$ with $i \neq 1,n$, we will assume that the Heegaard surface separates the two boundary components.  

Let $(\Sigma_0, H_0^-, H_0^+)$ be the Heegaard splitting formed by amalgamating the Heegaard splittings of $M_1$ and $M_2$, then amalagamating the resulting Heegaard splitting with the Heegaard splitting of $M_3$ and so on to $M_n$.  For each $i < n$, we will glue the boundary component of $M_i$ corresponding to $\partial_- H_i^-$ to the one boundary component of the manifold constructed up to that point.  For the last step, each of the two pieces has a single boundary component.  We will temporarily defy our convention and order the compression bodies in $M_n$ so that $H_n^+$ is the handlebody and $H_n^-$ is a compression body.  For each amalgamation, we will order the compression bodies so that they agree with the Heegaard splitting of the $M_i$ that we are attaching.

For $1 \leq i < n$, let $(\Sigma_i, H_i^-, H_i^+)$ be the Heegaard splitting resulting from the same construction except that at the step before we amalgamate with the Heegaard splitting for $M_{i+1}$, we will boundary stabilize the Heegaard splitting produced up to that point.

\begin{Lem}
\label{ex2boundlem}
The genus of $\Sigma_0$ is $3n + 1$, while the genus of each $\Sigma_i$ for $i \geq 1$ is $3n+2$.  For each $S_i \subset M_i$, the Heegaard genera of the components of $M \setminus S_i$ sum to $3n+4$.
\end{Lem}

\begin{proof}
The genera of $\Sigma_0,\dots,\Sigma_n$ can be calculated from the amalgamation constructions, as noted in Section~\ref{amalgsect}.  The Heegaard genus of an amalgamation is the sum of the Heegaard surface genera minus the sum of the genera of the amalgamating surfaces.  The Heegaard surface $\Sigma_0$ is an amalgamation of $n$ genus four Heegaard splittings along $n-1$ genus one surfaces, so its genus is $4n - 1(n-1) = 3n+1$.  Each $\Sigma_i$ for $i \geq 1$ is an amalagamation of $n-1$ genus four Heegaard surfaces and one genus five Heegaard surface along $n-1$ genus one surfaces, so each has Heegaard genus $4(n-1) + 5 - 1(n-1) = 3n+2$.

To calculate the Heegaard genera of the components of $M \setminus S_i$, consider a minimal genus Heegaard surface $\Sigma''$ for $M \setminus S_i$.  Compress the boundary as in Corollary~\ref{hgencoro} to a component of $\partial M_i$.  If $\Sigma'$ splits any $S_j$ then by Lemma~\ref{splitlem}, its genus is at least $3n+4$ and we're done.  Otherwise, $\Sigma'$ spans each $S_j$, so by the first half of Lemma~\ref{bothspanlem}, $\Sigma'$ is an amalgamation along each $S_j$.  By a calculation as above, the Heegaard genera of the two components sum to at least $3n+4$.  Conversely, one can always construct Heegaard surfaces of the components whose genera sum to $3n+4$.
\end{proof}

\begin{Lem}
For any $0 \leq i < j \leq n-1$, the stable genus of the Heegaard surfaces $\Sigma_i$ and $\Sigma_j$ is at least $3n+4$.
\end{Lem}

\begin{proof}
For each $i$, let $S_i \subset M$ be the image of the Heegaard surface for $M_i$ in $M$.  The image in $M$ of $\partial M_i$ is incompressible so $M_{S_i}$ is isotopic to the image of $M_i$.  By construction, $(\Sigma_0, H^-_0, H^+_0)$ spans each $(S_j, M_{S_j}^-, M_{S_j}^+)$ positively.  For $i > 0$, $(\Sigma_i, H^-_i, H^+_i)$ spans each $(S_j, M_{S_j}^-, M_{S_j}^+)$ positively for $j > i$ and negatively for $j \leq i$.  

Given any $k \neq l \leq n$, there will be some $i$ and $j$ such that $(S_k, M_{S_k}^-, M_{S_k}^+)$ and $(S_l, M_{S_l}^-, M_{S_l}^+)$ will span $(S_i, M_{S_i}^-, M_{S_i}^+)$ with the same sign, but span $(S_j, M_{S_j}^-, M_{S_j}^+)$ with opposite signs.  The same will be true of any common stabilization, so by Lemma~\ref{bothspanboundlem}, any common stabilization has genus at least $\min \{\frac{1}{2}(6n+8), 3n+4\} = 3n+4$.
\end{proof}

\begin{proof}[Proof of Theorem~\ref{mainthm2}]
By Lemma~\ref{ex2boundlem}, $\Sigma_0$ has genus $3n+1$ while for $0 < i \leq n$, each $\Sigma_i$ has genus $3n+2$.  For each $i \leq n$, let $\Sigma'_i$ be a genus $3n + 3$ stabilization of $\Sigma_i$.  If any of these two Heegaard surfaces are isotopic then they are genus $3n+3$ common stabilizations of some $\Sigma_k$, $\Sigma_j$.  Because $3n+3 < 3n+4$, this contradicts Lemma~\ref{ex2boundlem}.  We conclude that the surfaces $\Sigma_0,\dots,\Sigma_{n-1}$ are pairwise non-isotopic.
\end{proof}

\bibliographystyle{amsplain}
\bibliography{stabs2}

\providecommand{\bysame}{\leavevmode\hbox to3em{\hrulefill}\thinspace}
\providecommand{\MR}{\relax\ifhmode\unskip\space\fi MR }
\providecommand{\MRhref}[2]{%
  \href{http://www.ams.org/mathscinet-getitem?mr=#1}{#2}
}
\providecommand{\href}[2]{#2}
\begin{thebibliography}{10}

\bibitem{bachman}
D.~Bachman, \emph{{Stabilizations of Heegaard splittings of sufficiently
  complicated 3-manifolds (Preliminary Report)}}, preprint (2008),
  arXiv:0806.4689.

\bibitem{bsc:bndls}
D.~Bachman and S.~Schleimer, \emph{{Surface bundles versus Heegaard
  splittings}}, Communications in Analysis and Geometry 13 \textbf{5} (2005),
  1--26.

\bibitem{cass:red}
A.~Casson and C.~Gordon, \emph{{Reducing Heegaard splittings}}, Topology Appl.
  \textbf{27} (1987), 275--283.

\bibitem{cerf:strat}
J.~Cerf, \emph{La stratefacation naturelle des especes de fonctions
  differentiables reeles et la theoreme de la isotopie.}, Publ. Math. I.H.E.S.
  \textbf{39} (1970).

\bibitem{htt:stabs}
J.~Hass, A.~Thompson, and W.~Thurston, \emph{{Common stabilizations of Heegaard
  splittings}}, preprint (2008).

\bibitem{Hempel:complex}
J.~Hempel, \emph{3-manifolds as viewed from the curve complex}, Topology
  \textbf{40} (2001), no.~3, 631--657.

\bibitem{me:t3}
J.~Johnson, \emph{{Automorphisms of the 3-torus preserving a genus three
  Heegaard splitting}}, preprint (2007), arXiv:0708.2683.

\bibitem{inflects}
J.~Johnson, \emph{{Stable functions and common stabilizations of Heegaard
  splittings}}, preprint (2007), arXiv:0705.3712.

\bibitem{me:stabs}
J.~Johnson, \emph{{Flipping and stabilizing Heegaard splittings}}, preprint
  (2008), arXiv:0805.4422.

\bibitem{Kob:disc}
T.~Kobayashi and O.~Saeki, \emph{{The Rubinstein-Scharlemann graphic of a
  3-manifold as the discriminant set of a stable map.}}, Pacific Journal of
  Mathematics \textbf{195} (2000), no.~1, 101--156.

\bibitem{mather}
J.~Mather, \emph{{Stability of $C^\infty$ mappings V.}}, Advances in
  Mathematics \textbf{4} (1970), no.~3, 301--336.

\bibitem{morsedg}
Y.~Moriah and E.~Sedgwick, \emph{Heegaard splittings of twisted torus knots},
  preprint (2007).

\bibitem{reid}
K.~Reidemeister, \emph{{Zur dreidimensionalen Topologie}}, Abh. Math. Sem.
  Univ. Hamburg \textbf{11} (1933), 189--194.

\bibitem{rub:compar}
H.~Rubinstein and M.~Scharlemann, \emph{{Comparing Heegaard splittings of
  non-Haken 3-manifolds.}}, Topology \textbf{35} (1996), no.~4, 1005--1026.

\bibitem{sch:thin}
M.~Scharlemann and A.~Thompson, \emph{Thin position for 3-manifolds},
  Contemporary Mathematics \textbf{164} (1994), 231--238.

\bibitem{st:dist}
M.~Scharlemann and M.~Tomova, \emph{{Alternate Heegaard genus bounds
  distance}}, preprint (2005), ArXiv:math.GT/0501140.

\bibitem{sing}
J.~Singer, \emph{{Three-dimensional manifolds and their Heegaard diagrams}},
  Trans. Amer. Math. Soc. \textbf{35} (1933), 88--111.

\end{thebibliography}

\end{document}